\documentclass[number,seceqn,dvips]{arxbj}
\usepackage{cursive}  
\usepackage{graphicx}

\aid{0}
\volume{16}
\issue{2}
\pubyear{2010}
\firstpage{543}
\lastpage{560}
\doi{10.3150/09-BEJ222}

\makeatletter
\newtheorem{teorem}[lem]{Theorem}
\newremark{remark}{Remark}
\newtheorem{lem}{Lemma}[section]
\newremark{example}{Example}
\def\Th{\Theta}
\def\th{\theta}
\def\tth{\tilde{\theta}}
\def\hp{\hat{p}}
\def\I{\mathbf{I}}
\def\Z{\tilde{Z}}
\def\z{\tilde{z}}
\def\Y{\tilde{Y}}
\def\y{\tilde{y}}
\def\X{\tilde{X}}
\def\x{\tilde{x}}
\def\s2{v}
\newcommand{\C}{\vert }
\newcommand{\RR}{\mathbb{R}}
\newcommand{\EE}{E}
\makeatother

\begin{document}
\begin{frontmatter}

\title{Asymptotic minimax risk of predictive density estimation for non-parametric regression}
\runtitle{Asymptotic minimax density prediction}

\begin{aug}
\author[1]{\fnms{Xinyi} \snm{Xu}\thanksref{1}\ead[label=e1]{xinyi@stat.osu.edu}\corref{}} \and
\author[2]{\fnms{Feng} \snm{Liang}\thanksref{2}\ead[label=e2]{liangf@uiuc.edu}}
\runauthor{X. Xu and F. Liang}
\address[1]{Department of Statistics,  Ohio State University, 1958 Neil Avenue, Columbus, OH
43210-1247, USA. \mbox{\printead{e1}}}
\address[2]{Department of Statistics, University of Illinois at Urbana-Champaign, 725 S.
Wright Street, Champaign, IL 61820, USA. \printead{e2}}
\end{aug}

\received{\smonth{7} \syear{2008}}
\revised{\smonth{3} \syear{2009}}

\begin{abstract}
We consider the problem of estimating the predictive density of
future observations from a non-parametric regression model.
The~density estimators are evaluated under Kullback--Leibler divergence
and our focus is on establishing the exact asymptotics of minimax
risk in the case of Gaussian errors.  We derive the convergence rate
and constant for minimax risk among Bayesian predictive densities
under Gaussian priors and we show that this minimax risk is
asymptotically equivalent to that among all density estimators.
\end{abstract}


\end{frontmatter}

\section{Introduction}\label{sec1}

Consider the canonical non-parametric regression setup
\begin{eqnarray} \label{eq:X}
Y(t_i) = f(t_i) + \sigma \varepsilon_i,\qquad    i = 1, \ldots, n,
\end{eqnarray}
where $f$ is an unknown function in $\mathcal{L}^2[0, 1]$,
$t_i=i/n$ and the $\varepsilon_i$'s are i.i.d.~standard Gaussian random
variables.  We assume the noise level $\sigma$ is known and, without
loss of generality, set $\sigma = 1$ throughout.

Based on observing $Y=(Y(t_1), \dots, Y(t_n))$, estimating $f$ or
various functionals of $f$ has been the central problem in
non-parametric function estimation.  The~asymptotic optimality of
estimators is usually associated with the optimal rate of convergence in
terms of minimax risk.  A~huge body of literature has been devoted
to the evaluation of minimax risks under $\mathcal{L}^2$ loss over
certain function spaces; see, for example,  Pinsker \cite{r19},
Ibragimov and Has'minskii \cite{r14}, Golubev and Nussbaum
\cite{r12}, Efroimovich \cite{r6}, Belitser and Levit \cite{r2,r3}
and Goldenshluger and Tsybakov \cite{r11}. An excellent survey of the literature
in this area can be found in Efromovich \cite{r7}.

Sometimes, instead of estimating $f$ itself, one is interested in making
statistical inference about future observations from the same
process that generated $Y(t)$. A~predictive distribution function
assigns probabilities to all possible outcomes of a random variable.
It thus provides a complete description of the uncertainty
associated with a prediction.  The~minimaxity of predictive density
estimators has been studied for finite-dimensional parametric
models; see, for example, Liang and Barron \cite{r16},  George,
Liang and Xu \cite{r9}, Aslan \cite{r1} and George and Xu
\cite{r10}. However, so far, few results have been obtained on
predictive density estimation for non-parametric models.  The~major
thrust of this paper is to establish the asymptotic minimax risk for
predictive density estimation under Kullback--Leibler loss in the context of non-parametric regression.  Our result closely
parallels the well-known work by Pinsker \cite{r19} for
non-parametric function estimation under $\mathcal{L}^2$ loss and
provides a benchmark for studying the optimality of density
estimates for non-parametric regression.

Let $\Y=(\Y(u_1), \dots, \Y(u_m))^t$ denote a vector of future observations from model (\ref{eq:X}) at locations $\{u_j\}_{i=1}^m$.  To evaluate the
performance of density prediction across the whole curve, we assume that the $u_j$'s are equally spaced dense (that is, $m \ge n$) grids in $[0,1]$.  Given
$f$, the conditional density $p(\y \C f)$ is a product of $N(\y_j; f(u_j))$, where $N(\cdot; \mu)$ denotes a univariate Gaussian density function
with mean $\mu$ and unit variance.  Based on observing $Y=y$, we estimate $p(\y \C f)$ by a predictive density $\hp(\y \C y)$, a non-negative
function of $\y$ that integrates to 1 with respect to~$\y$.

Common approaches to constructing  $\hp(\y \C y)$ includes
the ``plug-in'' rule that simply substitutes an estimate $\hat{f}$ for $f$ in $p(\y \C f)$,
\begin{equation}  \label{hp:plug-in}
p(\y \C \hat{f}) = \prod_{j=1}^n N(\y_j; \hat{f}(u_j)),
\end{equation}
and the Bayes rule that integrates $f$  with respect to a prior $\pi$
to obtain
\begin{eqnarray} \label{hp:Bayes}
\int p(\y \C f) \pi(f \C y)\, \mathrm{d}f =  \frac{\int p(y \C f) p(\y \C f) \pi(f)\, \mathrm{d}f}{\int p(y \C f) \pi(f)\, \mathrm{d}f}.
\end{eqnarray}
We measure the discrepancy between $p(\y \C f)$ and $\hp(\y \C y)$
by the average Kullback--Leibler (KL) divergence
\begin{eqnarray}  \label{risk:simultaneous}
R(f, \hat{p}) =
\frac{1}{m} E_{Y, \tilde{Y} | f} \log \frac{p(\tilde{Y} \vert  f)}{\hp(\tilde{Y} \vert  Y)}.
\end{eqnarray}
Assuming that $f$ belongs to a function space $\mathcal F$, such as
a Sobolev space, we are interested in the minimax risk
\begin{equation} \label{minimax:f}
R(\mathcal{F}) = \min_{\hp} \max_{f \in \mathcal{F}}  R(f,
\hp).
\end{equation}

It is worth observing that in this framework, the densities of future
observations $(\Y_1, \dots, \Y_m)$ are estimated simultaneously by
$\hp(\y \C y)$. An alternative approach is to estimate the densities
individually by $\{\hp(\y_j \C y)\}_{j=1}^m$ with risk
\begin{equation} \label{risk:individual}
\frac{1}{m} \sum_{j=1}^m E_{Y, \Y | f}  \log \frac{p(\Y_j \C
f(u_j))}{\hp(\Y_j \C Y)}.
\end{equation}
When the $u_j$'s are equally spaced and $m$ goes to infinity, the
risk above converges to
\[
 \int_0^1 \EE_{Y, \Y | f} \log \frac{p(\Y \vert  f(u))}{\hp (\Y \C Y)} \,\mathrm{d} u,
 \]
which can be interpreted as the integrated KL risk of prediction at a random location $u$ in $[0,1]$.  This individual prediction problem can be
studied in our simultaneous prediction framework with $\hp(\y \C y)$ restricted to a product form, that is, $\hp(\y \C y) = \prod_{j=1}^m \hp(\y_j \vert
y)$.  For example, the plug-in estimator (\ref{hp:plug-in}) has such a product form and it is easy to check that its individual estimation risk
(\ref{risk:individual}) is the same as its simultaneous estimation risk (\ref{risk:simultaneous}). In general, simultaneous prediction considers a
broader class of $\hp$ than the one considered by individual prediction. Therefore, simultaneous prediction is more efficient since the corresponding
minimax risk (\ref{minimax:f}) is less than or equal to the one with individual prediction. This is distinct from estimating $f$ itself under
$\mathcal{L}^2$ loss where, due to the additivity of $\mathcal{L}^2$ loss, simultaneous estimation and individual estimation are equivalent.

This paper is organized as follows.  In Section \ref{sec2}, we show that the
problem of predictive density estimation for a non-parametric
regression model can be converted to the one for a Gaussian sequence
model with a constrained parameter space.  Direct evaluation of
the minimax risk  is difficult because of the
constraint on the parameter space.  Therefore, in Section \ref{sec3}, we first
derive the minimax risk over a special class of $\hp$ that consists
of predictive densities under Gaussian priors on the unconstrained
parameter space $\RR^n$.  Then, in
Section \ref{sec4}, we show that this minimax risk is asymptotically equivalent to the overall
minimax risk. Finally, in Section \ref{sec5}, we provide two explicit examples of minimax risks over $\mathcal{L}^2$ balls and Sobolev spaces.

\section{Connection to Gaussian sequence models}\label{sec2}

Let $\{\phi_i\}_{i=1}^\infty$ be the orthonormal trigonometric basis of
$\mathcal{L}^2[0, 1]$, that is,
\begin{eqnarray*}
\phi_0(t) \equiv 1,\qquad
\cases{
\phi_{2k-1} = \sqrt{2} \sin(2 \curpi k x),  \vspace*{2pt}\cr
\phi_{2k}   = \sqrt{2} \cos(2 \curpi k x),}
\qquad    k=1, 2, \dots.
\end{eqnarray*}
Then, $f = \sum_{i=1}^\infty \theta_i  \phi_i$,
where $\theta_i = \int_0^1 f(t) \phi_i(t) \,\mathrm{d}t$ is the coefficient
with respect to the $i$th basis element $\phi_i$. A~function space
$\mathcal F$ corresponds to a constraint on the parameter space of
$\theta$.  In this paper, we consider function spaces whose parameter spaces $\Th$ have ellipsoid constraints,
that is,
\begin{equation} \label{Theta(C)}
 \Th(C) = \Biggl\{ \th\dvtx  \sum_{i=1}^\infty  a_i^2 \th_i^2 \le C \Biggr\},
\end{equation}
where $a_1 \le a_2 \le \cdots$ and $a_n \rightarrow \infty$.

We approximate $f$ by a finite summation $f_n = \sum_{i=1}^n
\theta_i \phi_i$.  The~bias incurred by estimating $p(\y\vert f_n)$
instead of $p(\y\vert f)$ can be expressed as
\begin{eqnarray*}
\operatorname{Bias}(f, f_n) =  \frac{1}{m} \EE_{\Y |f}  \log \frac{p(\Y | f)}{p(\Y|f_n)}
 =   \frac{1}{2 m} \sum_{j=1}^m [f(u_j) - f_n(u_j)]^2
 =  \frac{1}{2m} \sum_{i=n+1}^{\infty} \theta_i^2.
\end{eqnarray*}
This bias is often negligible compared to the prediction risk (\ref{risk:simultaneous}); for example, it is of order $\mathrm{O}(n^{-2 \alpha})$ for Sobolov
ellipsoids $\Theta(C, \alpha)$, as defined in (\ref{Sobolev_space}).  Therefore, from now on, we set $f=f_n$.

Let $\theta = (\theta_1, \theta_2, \ldots, \theta_n)^t$, $\Phi_A$ be
a $n \times n$ matrix whose $(i,j)$th entry equals  $\phi_j(t_i)$
and $\Phi_B$ be a $m \times n$ matrix whose $(i,j)$th entry equals
 $\phi_j(u_i)$.  Then, $Y \C \theta$ and $\Y \C \theta$ are two
independent Gaussian vectors with $Y \C \theta \sim N (\Phi_A
\theta, \I_n)$ and $\Y \C \theta \sim N(\Phi_B \theta, \I_m)$, where $\I_n$ denotes the $n\times n$ identity matrix.  Note
that since the $t_i$'s and $u_j$'s are equally spaced, we have $\Phi_A^t \Phi_A = n
\I_n$ and $\Phi_B^t \Phi_B = m \I_n$.  Defining
\begin{equation} \label{eq:transform}
 X = \frac{1}{n} \Phi_A^t Y \quad  \mbox{and}\quad   \X = \frac{1}{m} \Phi_B^t \Y,
\end{equation}
it is then easy to check that $X$ and $\X$ are independent and
that\begin{equation} \label{model:normal}
 X \C \theta \sim N(\th, \s2_n \I_n)  \quad \mbox{and}\quad   \X \C \theta \sim N(\th, \s2_m \I_n),
\end{equation}
where $\s2_n=1/n$ and $\s2_m=1/m$.  We refer to the model
above as a \textit{Gaussian sequence model} since its number
of parameters is increasing at the same rate as the number of data
points.

Consider the problem of predictive density estimation for the Gaussian sequence model (\ref{model:normal}). Let $\hp(\x \C x)$ denote a predictive density function of $\x$ given $X=x$. The~incurred KL risk is defined to be
\[
R(\th, \hp) = \frac{1}{m} E_{X, \X \vert  \th} \log \frac{p(\X \C
\th)}{\hp(\X \C  X)}
\]
and the corresponding minimax risk is given by
\begin{equation} \label{def:R}
 R(\Theta) = \inf_{\hp} \sup_{\th \in \Th(C)} R(\th, \hp).
\end{equation}
The~following theorem states that the two minimax risks, the one associated with $(Y, \Y)$ from a non-parametric regression model and
the one associated with $(X, \X)$  from a normal sequence model, are equivalent.

\begin{teorem} \label{teo2.1}
$R({\mathcal F}) = R(\Theta)$, where $R({\mathcal F})$ is defined in \textup{(\ref{minimax:f})} and $R(\Theta)$ in \textup{(\ref{def:R})}.
\end{teorem}

\begin{pf}
See the \hyperref[append]{Appendix}.
\end{pf}

\begin{remark*}
The~idea of reducing a non-parametric
regression model to a Gaussian sequence model via an orthonormal
function basis has been widely used for non-parametric function
estimation.  Early references include Ibraginov and
Has'minskii \cite{r13}, Efromovich and Pinsker \cite{r8} and references therein. For recent developments, see Brown and Low \cite{r4},
Nussbaum \cite{r17,r18} and
Johnstone \cite{r15}.  Our proof of Theorem \ref{teo2.1}, given in the \hyperref[append]{Appendix}, implies that \textit{simultaneous}
estimation of predictive densities in these two models are
equivalent. However, this equivalence does not hold for the \textit{individual} estimation approach
described in Section \ref{sec1} because the product form of the density estimators, that is,
$\hp(\y \C y) = \prod_{j} \hp(\y_j \C y)$, is not retained under the transformation.
\end{remark*}

\section{\textit{Linear} minimax risk}\label{sec3}

Direct evaluation of the minimax risk (\ref{def:R}) is
difficult because the parameter space $\Theta(C)$ is constrained.  In this
section, we first consider a subclass of density estimators that have
simple forms and investigate the minimax risk over this subclass.
In next section, we then show that the minimax risk over this
subclass is asymptotically equivalent to the overall minimax risk~$R$.  Such an approach was first used in Pinsker
\cite{r19} to establish a minimax risk bound for the function
estimation problem. It inspired a series of developments, including
Belitser and Levit \cite{r2,r3}, Tsybakov \cite{r20} and
Goldenshluger and Tsybakov \cite{r11}.

Recall that in the problem of estimating the mean of a Gaussian
sequence model under $\mathcal{L}^2$ loss, diagonal linear
estimators of the form $\hat\theta_i = c_i x_i$ play an important
role.  Indeed, Pinsker \cite{r19} showed that when the parameter
space (\ref{Theta(C)}) is an ellipsoid, the minimax risk among
diagonal linear estimators is asymptotically minimax among all
estimators. Moreover, the results in Diaconis and Ylvisaker
 \cite{r5} imply that if such a diagonal linear estimator is  Bayes, then the prior~$\pi$ must be a Gaussian prior with a diagonal covariance matrix.  Similarly, in investigating the minimax risk of  predictive density estimation, we
first restrict our attention to a special class of~$\hp$ that
are Bayes rules under Gaussian priors over the
unconstrained parameter space $\RR^n$. Due to the above connection,
we call these predictive densities \textit{linear} predictive densities
and call the minimax risk over this class the \textit{linear}
minimax risk, even though `linear' does not have any literal
meaning in our setting.

Under a Gaussian prior $\pi_S(\theta)= N(0, S)$, where
$S=\operatorname{diag}(s_1, \dots, s_n)$ and $s_i \ge 0$ for $i = 1, \dots,
n$, the linear predictive density $\hp_S$ is given by
\begin{eqnarray} \label{eq:linear_estimator}
\hp_S(\x \vert  x) =  \int_{\RR^n} p(\x \vert  \theta) \pi_S(\theta \vert
x) \,\mathrm{d} \theta =  \frac{\int_{\RR^n} p(x \vert  \th) p(\x \vert  \th)
  \pi_S(\th) \,\mathrm{d} \th}{\int_{\RR^n} p(x \vert  \th) \pi_S(\th)\, \mathrm{d} \th}.
\end{eqnarray}
Note that $\hp_S$ is not a Bayes estimator for the problem described in Section \ref{sec2} because the prior
distribution $N(0,S)$ is supported on $\RR^n$ instead of on the
ellipsoidal space $\Theta$. Nonetheless, $\hp_S$ is a valid
predictive density function.

The~following lemma provides an explicit form of the average KL risk of $\hp_S$.

\begin{lem}\label{lem3.1}
The~average Kullback--Leibler risk (\ref{risk:simultaneous}) of $\hp_S$ is given by
\begin{equation} \label{eq:linear:risk}
R(\th, \hp_S) = \frac{n}{2 m} \log \frac{\s2_n}{\s2_{n+m}} +
\frac{1}{2m} \sum_{i=1}^n  \biggl[ \log \frac{\s2_{n+m} + s_i}{\s2_n +
s_i} + \frac{\s2_{n+m} + \th_i^2}{\s2_{n+m} + s_i} - \frac{\s2_n +
\th_i^2}{\s2_n + s_i}  \biggr],
\end{equation}
where $\s2_{n+m}=1/(n+m)$.
\end{lem}

\begin{pf}
Let $\hp_U$ denote the posterior predictive density under the
uniform prior $\pi_U \equiv 1$, namely,
\begin{eqnarray*}
\hp_U(\x \C x) =  \biggl( \frac{1}{2 \curpi \s2_{n + m}} \biggr)^{n/2}
\exp  \biggl(  -\frac{\|\x-x\|^2}{2 \s2_{n +m}}  \biggr).
\end{eqnarray*}
Then, by \cite{r9}, Lemma 2, the average KL
risk of $\hp_S$ is given by
\begin{eqnarray} \label{eq:marginal representation}
R(\th, \hp_S) = R(\th, \hp_U) - \frac{1}{m} E \log m_S(W; \s2_{n+m})
+ \frac{1}{m} E \log m_S(X; \s2_n),
\end{eqnarray}
where
\[
W  = \frac{\s2_m X + \s2_n \tilde{X}}{\s2_{n+m}} \sim N(\theta, \s2_{n+m} I)
\]
and $m_S (x; \sigma^2)$ denotes the marginal distribution of $X \C
\theta \sim N_n(\theta, \sigma^2I)$ under the normal prior~$\pi_S$.
It is easy to check that
\begin{eqnarray} \label{lemma1:part1}
R(\th, \hp_U) &=& \frac{1}{m} \EE \log \frac{p(\x \C
\theta)}{\hp_U(\x | x)} = \frac{n}{2 m} \log
\frac{\s2_n}{\s2_{n+m}}
\end{eqnarray}
and
\begin{eqnarray}
E \log m_S(W; \s2_{n+m}) &=& -\frac{n}{2m} \sum_{i=1}^n \log [ 2 \curpi (\s2_{n+m} + s_i ) ]
- \frac{1}{2m} \sum_{i=1}^n \frac{\s2_{n+m} + \th_i^2}{\s2_{n+m} + s_i}, \label{lemma1:part2} \\
E \log m_S(X; \s2_n) &=& -\frac{n}{2m} \sum_{i=1}^n \log [ 2
\curpi (\s2_n + s_i ) ] - \frac{1}{2m} \sum_{i=1}^n
\frac{\s2_n + \th_i^2}{\s2_n + s_i}. \label{lemma1:part3}
\end{eqnarray}
The~lemma then follows immediately by combining
equations (\ref{eq:marginal representation})--(\ref{lemma1:part3}).
\end{pf}

We denote the linear minimax risk over all $\hp_S$ by $R_L(\Theta)$,
that is,
\begin{equation} \label{R_L:def}
 R_L(\Theta) =  \inf_S \sup_{\th \in \Th(C)} R(\th, \hp_S).
\end{equation}
This linear minimax risk is not directly tractable because the
inside maximization is over a constrained space $\Theta(C)$.  In the following theorem, we first show that we can switch the order
of $\inf$ and $\sup$ in equation (\ref{R_L:def}) and then evaluate
$R_L$ using the Lagrange multiplier method.

The~following notation will be useful throughout.  Let
$\tilde{\lambda}(C, \s2_n, \s2_{n+m})$ denote a solution of the
equation
\begin{equation} \label{def:tl}
 \sum_{i=1}^n a_i^2  \Biggl[  (\s2_n-\s2_{n+m}  ) \sqrt{1 + \frac{4 \tilde{\lambda}/a_i^2}{
 \s2_n-\s2_{n+m}}} -  (\s2_n+\s2_{n+m}  ) \Biggr]_+ = 2 C,
\end{equation}
where $[ x ]_{+} = \sup(x, 0)$, and let $\tth_i^2$ be
\begin{equation}  \label{def:tth}
 \tth_i^2 = \frac{1}{2}  \Biggl[ (\s2_n-\s2_{n+m}  ) \sqrt{1
 +\frac{4 \tilde{\lambda}/a_i^2}{\s2_n-\s2_{n+m}}} -  (\s2_n+\s2_{n+m}
  ) \Biggr]_+
\end{equation}
for $i=1, 2, \dots, n$.

\begin{teorem} \label{teo3.2}
Suppose that the parameter space $\Theta(C)$ is an ellipsoid, as defined
in (\ref{Theta(C)}).  The~linear minimax risk is then given by
\begin{eqnarray}
R_L(\Th)  &=&  \inf_S \sup_{\th \in \Th(C)} R(\th, \hp_S) =
\sup_{\th \in \Th(C)} \inf_S R(\th, \hp_S) \label{eq:exchange} \\
&=&\frac{n}{2m} \log \frac{\s2_n}{\s2_{n+m}} + \frac{1}{2m}
\sum_{i=1}^n \log \frac{ \s2_{n+m} + \tth_i^2}{\s2_n + \tth_i^2},
\label{eq:solution}
\end{eqnarray}
where $\tth_i^2$ is defined as in (\ref{def:tth}).  The~linear
minimax estimator $\hp_{\tilde{V}}$ is the Bayes predictive density
under a Gaussian prior
\begin{equation}
\pi_{\tilde{V}} (\theta) = N(0, \tilde{V}),\qquad  \mbox{ where } \tilde{V}
= \operatorname{diag}(\tth_1^2, \tth_2^2, \ldots, \tth_n^2),  \label{pi_V}
\end{equation}
namely,
\begin{eqnarray*}
\hp_{\tilde{V}}(\x \C x) = N(\theta_{\tilde{V}}, \Sigma_{\tilde{V}}),
\end{eqnarray*}
with
\begin{eqnarray*}
\theta_{\tilde{V}} &=& \biggl(\frac{\tth_1^2}{\tth_1^2 + v_n} x_1, \ldots, \frac{\tth_n^2}{\tth_n^2 + v_n} x_n\biggr)', \\
\Sigma_{\tilde{V}} &=& \operatorname{diag}\biggl(\frac{\tth_1^2 v_n}{\tth_1^2 + v_n}
+ v_m, \ldots, \frac{\tth_n^2 v_n}{\tth_n^2 + v_n} + v_m \biggr).
\end{eqnarray*}
\end{teorem}

\begin{pf}
We first prove equality (\ref{eq:solution}). It is easy to
check that for any fixed $\theta$, $R(\theta, \hp_S)$ achieves its
minimum at $S=\operatorname{diag}(\theta_1^2, \dots, \theta_n^2)$, and
\[
\inf_S R(\th, \hp_S) =
\frac{n}{2m} \log \frac{\s2_{n+m}}{\s2_n} + \frac{1}{2m}
\sum_{i=1}^n \log \frac{ \s2_{n+m} + \th_i^2}{\s2_n + \th_i^2}.
\]
To calculate the maximum of the above quantity  over $\theta \in
\Theta(C)$, one needs to solve
\[
\sup  \Biggl\{ \sum_{i=1}^n \log \frac{
  \s2_{n+m} + \th_i^2}{\s2_n + \th_i^2}  \dvtx   \sum_{i=1}^n a_i^2 \th_i^2 \le C
   \Biggr\}.
    \]
With the Lagrangian
\[ \mathcal{L} =\sum_{i=1}^n \log \frac{
  \s2_{n+m} + \th_i^2}{\s2_n + \th_i^2} - \frac{1}{\lambda}  \Biggl( \sum_{i=1}^n
  a_i^2 \th_i^2 - C \Biggr),
\]
simple calculation reveals that the maximum is attained at $\tth_i$
given by (\ref{def:tth}).

Next, we prove equality (\ref{eq:exchange}), that is, that the order of inf
and sup can be exchanged.  Note that for any diagonal matrix
$\tilde{S}$, we have
\begin{equation} \label{ineq:chain}
\sup_{\theta \in \Theta(C)} R(\hat
 p_{\tilde{S}}, \theta) \ge \inf_S \sup_{\theta \in
\Theta(C)} R(\hat p_S, \theta) \ge \sup_{\theta \in \Theta(C)}
\inf_S R(\hat p_S, \theta).
\end{equation}
Therefore, if there exists an $\tilde{S}$ such that
\[
\sup_{\theta \in \Theta(C)} R(\hat p_{\tilde{S}}, \theta) - \sup_{\theta \in \Theta(C)} \inf_S R(\hat p_S, \theta) \leq 0 ,
\]
then all of the inequalities in (\ref{ineq:chain}) become equalities.

If we let $\tilde{S} = \operatorname{diag}(\tth_1^2, \ldots, \tth_n^2)$, then
\begin{eqnarray*}
R(\hat p_{\tilde{S}}, \theta)- \sup_{\theta \in \Theta(C)} \inf_S
R(\hp_S, \theta) & =&  \frac{1}{2m} \sum_{i=1}^n \frac{(\s2_n -
\s2_{n+m})(\theta_i^2 - \tilde{\theta}_i^2)}{(\s2_n +
\tilde{\theta}_i^2)(\s2_{n+m} + \tilde{\theta}_i^2)}\\
&=&  \frac{1}{2m} \frac{\sum_{i=1}^n  a_i^2 \theta_i^2
-C}{\tilde{\lambda}},
\end{eqnarray*}
where the second equality holds because $\sum_{i=1}^n a_i^2 \tth_i^2
= C$ and $\tilde{\theta}_i^2$ is a solution to
\begin{eqnarray*}
\frac{\partial \mathcal{L}}{\partial \theta_i^2} = \frac{\s2_n -
  \s2_{n+m}}{(\s2_n + \theta_i^2)(\s2_{n+m} + \theta_i^2)} - \frac{a_i^2}{\tilde{\lambda}} =
0.
\end{eqnarray*}
Since $\theta \in \Theta(C)$ implies that $\sum_{i=1}^n a_i^2
\th_i^2  \le C$,  we have
\[
\sup_{\theta \in \Theta(C)} R(\hat p_{\tilde{S}}, \theta) -
\sup_{\theta \in \Theta(C)} \inf_S R(\hp_S, \theta) \le
 \frac{1}{2m} \frac{C-C}{\tilde{\lambda}} = 0,
 \]
which completes the proof.
\end{pf}

\begin{remark*}  Note that $a_1 \leq a_2 \leq \cdots$, so we
have $\tth_i^2=0$ for $i > N$, where
\begin{equation} \label{def:N}
 N = \sup   \biggl\{i\dvtx  a_i^2 \le  \tilde{\lambda}  \biggl( \frac{1}{\s2_{m+n}}-\frac{1}{\s2_n}  \biggr) =
m \tilde{\lambda}  \biggr\}.
\end{equation}
This implies that the prior distribution corresponding to the linear
minimax estimator, that is, $\pi_{\tilde{V}}(\theta) = \prod_{i=1}^n
N(0, \tilde{\theta}_i^2)$, puts a point mass at zero for $\theta_i$
for all $i >N$.
\end{remark*}


\section{Asymptotic minimax risk}\label{sec4}

In this section, we turn to establishing the asymptotic behavior of
the minimax risk $R(\Th)$ over all predictive density estimators. By
definition, $R(\Th) \le R_L(\Th)$.  We extend the approach in
 \cite{r2} to show that the difference between
$R(\Th)$ and $R_L(\Th)$ vanishes as the number of observations~$n$
goes to infinity.  Therefore, the overall minimax risk is
asymptotically equivalent to the linear minimax risk.  This also
implies that the Gaussian prior $\pi_{\tilde{V}}$ defined in
(\ref{pi_V}) is asymptotically least favorable.

The~following lemma provides a lower bound for the overall minimax risk $R(\Th)$ under some conditions.

\begin{lem} \label{lem4.1}
Let $\{s_i^2\}_{i=1}^n$ be a sequence such that for some $\alpha >
0$,
\begin{equation} \label{cond:s_i}
\sum_{i=1}^n a_i^2 s_i^2 +  \Biggl[ - 8 \alpha  \Biggl(\sum_{i=1}^n
a_i^4 s_i^4  \Biggr) \log \s2_n  \Biggr]^{1/2} \le C.
\end{equation}
Then, as $n \rightarrow \infty$, the minimax risk $R(\Theta)$ has the
following lower bound:
\[ R(\Th) \ge  \frac{n}{2m} \log \frac{\s2_n}{\s2_{n+m}} + \frac{1}{2m}\sum_{i=1}^n \log \frac{\s2_{n+m} +
  s_i^2}{\s2_n + s_i^2} + \mathrm{O}(\s2_n^{\alpha}). \]
\end{lem}

\begin{pf}
See the \hyperref[append]{Appendix}.
\end{pf}

Note that, as shown in the proof, for a posterior density with a
Gaussian prior $\pi_S = N(0,S)$, where $S=\operatorname{diag}(s_1, \dots, s_n)$,
condition (\ref{cond:s_i}) guarantees $\pi_S$ to have most of its
mass inside $\Th$, in the sense that $\pi_S(\Th^c) \le v_n^{2
\alpha}$ for some $\alpha > 0$.

With the lower bound in the above lemma, we are ready to prove
the main result in this paper, which shows that the overall minimax
risk $R(\Theta)$ is asymptotically equivalent to the linear minimax
risk $R_L(\Theta)$.

\begin{teorem} \label{teo4.2}
Suppose that $\Theta$ is the ellipsoid defined in (\ref{Theta(C)}) and
$\tth^2$ is defined in (\ref{def:tth}).  If $m = \mathrm{O}(n)$ and
\begin{equation} \label{eq:alpha_v}
\log(1/\s2_n) \sum_{i=1}^n a_i^4 \tilde{\theta}_i^4 = \mathrm{o}(1),\qquad
\mbox{as } \s2_n \rightarrow 0,
\end{equation}
then
\begin{eqnarray}
\lim_{\s2_n \rightarrow 0} \frac{R(\Theta)}{R_L(\Theta)} = 1.
\end{eqnarray}
\end{teorem}

\begin{pf}
By definition, $R(\Theta) \leq R_L(\Theta)$. So, to
prove this theorem, it suffices to show that as \mbox{$\s2_n \rightarrow
0$},
\[
R(\Theta) \geq R_L(\Theta)\bigl(1 - \mathrm{o}(1)\bigr).
\]

For a fixed constant $\alpha > 1$, let $\gamma =  \frac{1}{C} [
8 \alpha \log(1/\s2_n) \sum_{i=1}^n
  a_i^4 \tilde{\theta}_i^4 ]^{{1}/{2}}$ and let $b_i^2 = \tilde{\theta}_i^2 (1 + \gamma)^{-1}$ for $i=1, \dots, n$.   It is easy to check that the sequence $\{b_i\}_{i=1}^n$
  satisfies the condition (\ref{cond:s_i}). Therefore, by Theorem  \ref{lem4.1},
\begin{eqnarray}\label{ineq:lower_bound}
R(\Theta) &\geq& \frac{n}{2m} \log \frac{\s2_n}{\s2_{n+m}} + \frac{1}{2m} \sum_{i=1}^n
\log \frac{\s2_{n+m} + b_i^2}{\s2_n + b_i^2} + \mathrm{O}(\s2_n^\alpha) \nonumber \\[-8pt]\\[-8pt]
            & =  & R_L(\Theta) - \frac{1}{2m} \sum_{i=1}^n \log
            \frac{(\s2_n + b_i^2)(\s2_{n+m} + \tilde{\theta}_i^2)}{(\s2_{n+m} +
            b_i^2)(\s2_n + \tilde{\theta}_i^2)} + \mathrm{O}(\s2_n^\alpha)\qquad
\mbox{as } \s2_n \rightarrow 0.\nonumber
\end{eqnarray}
Next, we will derive the convergence rate of $R_L(\Theta)$ and
show that the other terms are of smaller order.

Using the fact that $\tilde{\th}_i^2=0$ for $i > N$ (see
(\ref{def:N})), we can rewrite $R_L(\Theta)$ as
\begin{eqnarray*}
R_L(\Theta) &=& \frac{n}{2m} \log \frac{\s2_n}{\s2_{n+m}} + \frac{1}{2m} \sum_{i=1}^N \log \frac{\s2_{n+m} + \tilde{\theta}_i^2}{\s2_n + \tilde{\theta}_i^2} + \frac{1}{2m} \sum_{i = N}^n \log \frac{\s2_{n+m}}{\s2_n} \\
            &=& \frac{1}{2m} \sum_{i=1}^N \log \frac{(\s2_{n+m} + \tilde{\theta}_i^2) \s2_n}{(\s2_n + \tilde{\theta}_i^2) \s2_{n+m}} \\
            &=& \frac{1}{2m} \sum_{i=1}^N \log \biggl( 1 + \frac{(\s2_n - \s2_{n+m}) \tilde{\theta}_i^2}{(\s2_n + \tilde{\theta}_i^2) \s2_{n+m}} \biggr).
\end{eqnarray*}
When $m = \mathrm{O}(n)$, we have $\s2_n - \s2_{n+m} = \mathrm{O}(\s2_n)$ and $\s2_n +
\s2_{n+m} = \mathrm{O}(\s2_n)$.  Therefore, by means of a Taylor expansion,
\begin{eqnarray} \label{eq:R_L_order}
R_L = \mathrm{O}\Biggl( \frac{1}{2m} \sum_{i=1}^N \frac{(\s2_n - \s2_{n+m})
\tilde{\theta}_i^2}{(\s2_n + \tilde{\theta}_i^2) \s2_{n+m}} \Biggr)
\ge \mathrm{O}\biggl(\frac{1}{m}\biggr).
\end{eqnarray}

Similarly, since $b_i = \tilde{\th}_i^2 = 0$ for $i > N$, the second
term in (\ref{ineq:lower_bound}) can be written as
\begin{eqnarray*}
\frac{1}{2m} \sum_{i=1}^n \log \frac{(\s2_n + b_i^2)(\s2_{n+m} +
\tilde{\theta}_i^2)}{(\s2_{n+m} + b_i^2)(\s2_n +
\tilde{\theta}_i^2)} = \frac{1}{2m} \sum_{i=1}^N \log \frac{(\s2_n +
b_i^2)(\s2_{n+m} + \tilde{\theta}_i^2)}{(\s2_{n+m} + b_i^2)(\s2_n +
\tilde{\theta}_i^2)}.
\end{eqnarray*}
For every $1 \le i \le N$, we have
\begin{eqnarray*}
\log \frac{(\s2_n + b_i^2)(\s2_{n+m} +
\tilde{\theta}_i^2)}{(\s2_{n+m} + b_i^2)(\s2_n
  + \tilde{\theta}_i^2)} &  = & \log \biggl( \frac{[ (1 +
    \gamma) \s2_n + \tilde{\theta}_i^2 ] (\s2_{n+m} +
  \tilde{\theta}_i^2)}{[ (1 + \gamma) \s2_{n+m} +
    \tilde{\theta}_i^2 ] (\s2_n + \tilde{\theta}_i^2)} \biggr)\\
        &  = & \log \biggl( 1 + \gamma \frac{(\s2_n - \s2_{n+m})
  \tilde{\theta}_i^2}{(\s2_n + \tilde{\theta}_i^2)(\s2_{n+m} +
  \tilde{\theta}_i^2) + \gamma \s2_n (\s2_{n+m} + \tilde{\theta}_i^2)}
\biggr) \\
& \leq & \log \biggl( 1 + \gamma \frac{(\s2_n - \s2_{n+m})
  \tilde{\theta}_i^2}{(\s2_n + \tilde{\theta}_i^2)\s2_{n+m}}
\biggr).
\end{eqnarray*}
Again using a Taylor expansion, as well as the condition that
$\gamma=\mathrm{o}(1)$, we obtain
\begin{eqnarray} \label{eq:term_2}
\frac{1}{2m} \sum_{i=1}^n \log \biggl( 1 + \gamma \frac{(\s2_n -
\s2_{n+m}) \tilde{\theta}_i^2}{(\s2_n +
\tilde{\theta}_i^2)\s2_{n+m}} \biggr) &  = & \mathrm{O}\Biggl(
\frac{\gamma}{2m} \sum_{i=1}^n  \frac{(\s2_n - \s2_{n+m})
\tilde{\theta}_i^2}{(\s2_n + \tilde{\theta}_i^2)\s2_{n+m}}\Biggr) =
\mathrm{o}(R_L).
\end{eqnarray}

Finally, since $m=\mathrm{O}(n)$, by choosing $\alpha > 1$, the last term in
(\ref{ineq:lower_bound}) satisfies
\begin{eqnarray} \label{eq:term_3}
\s2_n^\alpha = \mathrm{o}(1).
\end{eqnarray}

Combining (\ref{ineq:lower_bound})--(\ref{eq:term_3}), the theorem then follows.
\end{pf}

\section{Examples}\label{sec5}

In this section, we apply Theorems \ref{teo3.2} and \ref{teo4.2} to
establish asymptotic behaviors of minimax risks over some
constrained parameter spaces.
In particular, we consider the asymptotics over $\mathcal{L}^2$ balls and Sobolev ellipsoids.

\begin{example}\label{ex1}
Suppose that $m=n$ and $\theta$ is restricted
in an $\mathcal{L}^2$ ball,
\begin{eqnarray} \label{L2_ball}
\Theta(C) = \Biggl\{\theta\dvtx  \sum_{i=1}^n \theta_i^2 \leq C\Biggr\}.
\end{eqnarray}
The~$\mathcal{L}^2$ ball can be considered as a variant of the ellipsoid (\ref{Theta(C)}) with $a_1 = a_2 = \cdots = a_n = 1$ and $a_{n+1} =
a_{n+2} = \cdots = \infty$.  Although the values of the $a_i$'s here depend on $n$, the proofs of the above theorems are still valid.  It is easy to
see that $N$ defined in (\ref{def:N}) is equal to $n$ and that $\tth_1^2 = \tth_2^2 = \cdots = \tth_n^2 = C/n$. Therefore,
\begin{eqnarray*}
(\log n) \sum_{i=1}^n a_i^4 \tth_i^4 = (\log n) \cdot \frac{C^2}{n}
= \mathrm{o}(1).
\end{eqnarray*}
By Theorem \ref{teo4.2}, the minimax risk among all predictive density
estimators is asymptotically equivalent to the minimax risk among
\textit{linear} density estimators.  Furthermore, by Theorem \ref{teo3.2},
\begin{eqnarray}
\lim_{n \rightarrow \infty} R(\Theta(C)) = \lim_{n \rightarrow \infty} R_L(\Theta(C)) = \frac{1}{2} \log 2 + \frac{1}{2} \log \frac{{1}/{(2n)} +
{C}/{n}}{{1}/{n} + {C}/{n}} \nonumber = \frac{1}{2} \log \frac{1+2C}{1+C}.
\end{eqnarray}

Note that this minimax risk is strictly smaller than the minimax risk over the class of plug-in estimators since, for any plug-in density $\hp(\x \C
\hat\theta)$,
\begin{eqnarray} \label{plugin_risk}
R(\theta, \hp) = \frac{1}{n} E \log \frac{p(\x \C \theta)}{p(\x \C \hat\theta)} = \frac{1}{n} E \biggl[- \frac{\|x-\theta\|^2 -
\|x-\hat\theta\|^2}{2/n} \biggr] = \frac{1}{2} E \|\hat\theta - \theta\|^2
\end{eqnarray}
and by Pinsker's theorem, the minimax risk of estimating $\theta$ under squared error loss is $C/(1+C)$, which is larger than
$ \log \frac{1+2C}{1+C}$, by the fact that $x > \log(1+x)$ for any $x>0$.
\end{example}

\begin{example}\label{ex2}
Suppose that $m=n$ and $\theta$ is restricted
in a Sobolev ellipsoid
\begin{eqnarray} \label{Sobolev_space}
\Theta(C, \alpha) = \Biggl\{\theta\dvtx  \sum_{i=1}^\infty a_i^2 \theta_i^2 \leq C\Biggr\},
\end{eqnarray}
where $a_{2i} = a_{2i-1} = (2i)^\alpha (\alpha
> 0)$ for $i=1, 2, \dots.$  Then, by (\ref{def:N}), we have $a_N^2 / \tilde{\lambda} n \sim N^{2 \alpha} / \tilde{\lambda} n \rightarrow 1$ as $n \rightarrow \infty$.  Substituting this relation into equation
(\ref{def:tl}) yields
\begin{eqnarray*}
2C &\sim& \sum_{i=1}^N i^{2 \alpha} \biggl(\frac{1}{2n} \sqrt{1+8 \tilde{\lambda} n i^{-2 \alpha}} - \frac{3}{2n} \biggr) \\
   &=& \frac{1}{2n} \sum_{i=1}^N i^{2 \alpha} \bigl(\sqrt{1+8 N^{2\alpha} i^{-2 \alpha}} - 3 \bigr)\bigl(1 + \mathrm{o}(1)\bigr).
\end{eqnarray*}
Using the Taylor expression
\begin{eqnarray*}
\sqrt{1+8 N^{2\alpha} i^{-2 \alpha}} = \sum_{k=0}^\infty \frac{2
\sqrt{2} (-1)^k (2k)!}{(1-2k) k!^2 32^k} \biggl( \frac{i}{N}
\biggr)^{(2k-1) \alpha}
\end{eqnarray*}
and the asymptotic relation
\begin{eqnarray*}
\sum_{i=1}^N i^r &=& \frac{N^{r+1}}{r+1}\bigl(1+\mathrm{o}(1)\bigr)\qquad    \mbox{as }
N \rightarrow \infty,  r > -1,
\end{eqnarray*}
we obtain
\begin{eqnarray*}
N = M n^{{1}/{(2\alpha+1)}}\bigl(1+\mathrm{o}(1)\bigr)\quad    \mbox{and}\quad
\tilde{\lambda} = M n^{{- 2\alpha}/{(2\alpha+1)}}\bigl(1+\mathrm{o}(1)\bigr),
\end{eqnarray*}
where
\[
M = \Biggl[ 4C\Big/ \Biggl( \sum_{k=0}^\infty \frac{2
\sqrt{2} (-1)^k (2k)!}{(1-2k) k!^2 32^k } \cdot
\frac{1}{(2k+1)\alpha + 1}  - \frac{3}{2 \alpha+1} \Biggr)
\Biggr]^{{1}/{(2\alpha+1)}}.
\]
  Note that, by (\ref{def:tth}),
\begin{eqnarray*}
\tilde{\theta}_i^2 = \frac{1}{2} \Biggl[
\frac{1}{2n}\sqrt{1+8\biggl(\frac{N}{i}\biggr)^{2 \alpha}} -
\frac{3}{2n}\Biggr]_+ \bigl(1+\mathrm{o}(1)\bigr).
\end{eqnarray*}
Therefore,
\begin{eqnarray*}
(\log n) \sum_{i=1}^N a_i^4 \tilde{\theta}_i^4  = \mathrm{O}\biggl((\log n)
\cdot \frac{N^{4\alpha+1}}{n^2}\biggr) = \mathrm{O}\bigl((\log n) \cdot
n^{-{1}/{(2\alpha+1)}}\bigr) = \mathrm{o}(1).
\end{eqnarray*}
By Theorem \ref{teo4.2}, the minimax risk among all predictive density
estimators is asymptotically equivalent to the minimax risk among
the \textit{linear} density estimators.  Furthermore, by Theorem~\ref{teo3.2},
\begin{eqnarray*}
R_L(\Theta(C, \alpha))  &=& \frac{1}{2} \log 2 + \frac{1}{2n} \sum_{i=1}^N \log \frac{{1}/{(2n)}+\tilde{\theta}_i^2}{{1}/{n}+\tilde{\theta}_i^2} +
\frac{n-N}{2n}\log \frac{1}{2} \\
                 &=& \frac{1}{2n} \sum_{i=1}^N \log \frac{{1}/{n}+2\tilde{\theta}_i^2}{{1}/{n}+\tilde{\theta}_i^2} \\
                 &=& \frac{1}{2n} \sum_{i=1}^N \log \biggl(1 + \frac{\tilde{\theta}_i^2}{{1}/{n}+\tilde{\theta}_i^2} \biggr).
\end{eqnarray*}
It is difficult to calculate an explicit form of the optimal
constant for the minimax risk due to the $\log$ function, but we can
get an accurate bound for it. By Taylor expansion, there exists
$x_i^* \in (0, 1), i = 1, 2, \dots, N$, such that
\begin{eqnarray*}
R_L(\Theta(C, \alpha))  =  \frac{1}{2n} \sum_{i=1}^N
\biggl(\frac{1}{1+x_i^*}
\frac{\tilde{\theta}_i^2}{{1}/{n}+\tilde{\theta}_i^2} \biggr)
\in \Biggl(\frac{1}{4n}\sum_{i=1}^N
\frac{\tilde{\theta}_i^2}{{1}/{n}+\tilde{\theta}_i^2},
\frac{1}{2n}\sum_{i=1}^N
\frac{\tilde{\theta}_i^2}{{1}/{n}+\tilde{\theta}_i^2}
                       \Biggr).
\end{eqnarray*}
Moreover,
\begin{eqnarray*}
\sum_{i=1}^N
\frac{\tilde{\theta}_i^2}{{1}/{n}+\tilde{\theta}_i^2} =
\sum_{i=1}^N \frac{{1}/{(4n)} \sqrt{1+8({N}/{i})^{2 \alpha}} -
{3}/{(4n)}}{{1}/{n} + {1}/{(4n)} \sqrt{1+8({N}/{i})^{2
\alpha}} - {3}/{(4n)}} = K \cdot N,
\end{eqnarray*}
where
\[
K = 1 +
 \frac{1}{2(2\alpha + 1)} - \frac{1}{2} \sum_{k=0}^\infty \frac{2 \sqrt{2} (-1)^k (2k)!}{(1-2k) k!^2 32^k } \cdot
\frac{1}{(2k+1)\alpha + 1} .
\]
  Therefore,
\begin{eqnarray}
\lim_{n\rightarrow \infty} n^{{2 \alpha}/{(2 \alpha + 1)}}
R(\Theta(C, \alpha)) = \lim_{n \rightarrow \infty} n^{{2
\alpha}/{(2 \alpha + 1)}} R_L(\Theta(C, \alpha)) \in \bigl(\tfrac{1}{4}
KM, \tfrac{1}{2} KM \bigr),
\end{eqnarray}
that is, the convergence rate is $n^{-{2 \alpha}/{(2 \alpha + 1)}}$ and the convergence constant is between $\frac{1}{4}KM$ and~$\frac{1}{2}KM$.
\end{example}

\begin{figure}

\includegraphics{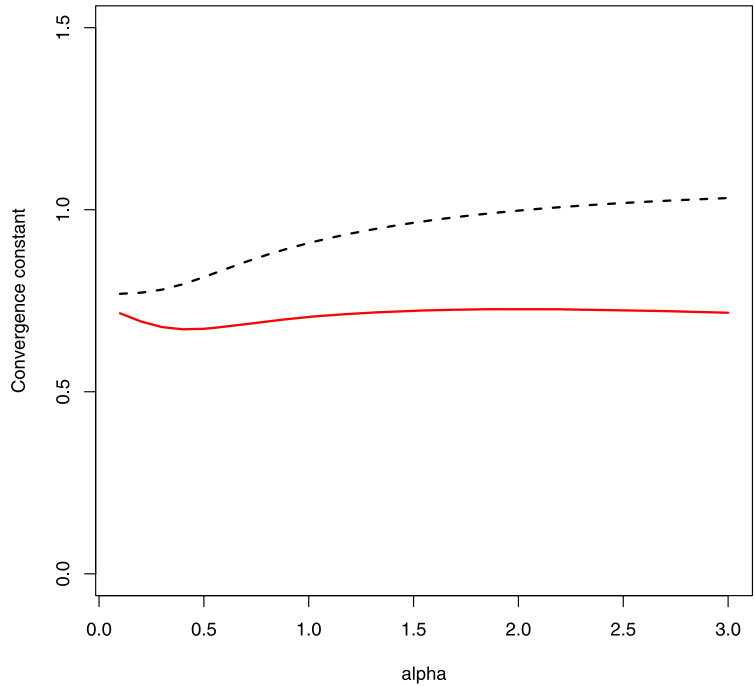}

\caption{Convergence constants of the overall minimax risk (the lower red line)
and the minimax risk over the class of plug-in estimators (the upper red line).
Here, the sample size $n=10\,000\,000$ and $C=1$. }\label{fig1}
\end{figure}

As in Example \ref{ex1}, we compare the asymptotics of this minimax risk with the one over the class of plug-in estimators, where the latter can be easily
computed by (\ref{plugin_risk}) and the results in  \cite{r19}.  Direct comparison reveals that the convergence rates of both minimax risks are
$n^{2\alpha/(2\alpha + 1)}$ and the convergence constants can both be written in the form
$C^{1/(2\alpha+1)} f(\alpha)$, where $f(\alpha)$ is a
function depending only on $\alpha$.  Although it is hard to obtain an explicit representation for the convergence constant for the overall minimax risk,
our simulation result in Figure \ref{fig1} shows that it is strictly smaller than that over the class of plug-in estimators.

\begin{appendix}\label{append}

\section*{Appendix: Proofs}

\renewcommand{\theequation}{A.\arabic{equation}}
\setcounter{equation}{0}
In this appendix, we provide the proofs of Theorem \ref{teo2.1} and Lemma \ref{lem4.1}.

\begin{pf*}{Proof of Theorem \ref{teo2.1}}  Let $\Psi$ be an
$m \times m$ matrix whose $(i, j)$th entry equals
$\phi_{j}(u_i)$. Since the $\phi_j$'s form an
orthogonal basis for $\mathcal{L}^2$ and the $u_i$'s are equally spaced,
we have $\Psi^t \Psi = \I_{m}$. Consider the transformation $\frac{1}{m}\Psi^t \Y$. Since the first $n$ columns of $\Psi$ are
$\Phi_B$, the first $n$ elements of the transformed vector are just $\X$, defined in (\ref{eq:transform}),
and we denote the remaining $(m-n)$ elements by $\Z$.  It is easy to check that $\X \C \theta
\sim N_n(\theta,   \frac{1}{m} \I_n)$ and $\Z \sim N_{m-n}(\mathbf{0},   \frac{1}{m} \I_{m-n})$
are independent multivariate Gaussian variables, and the target
density function $p(\y \C f)$ satisfies
\begin{eqnarray}
p(\y \C f) = p(\x, \z \C \theta) J_{\x, \z}(\y), \label{trans_1}
\end{eqnarray}
where $J_{\x, \z}(\y)$ is the Jacobian for this transformation.
Similarly, any predictor density estimator $\hp(\y \C y)$ can be
rewritten as
\begin{eqnarray}
\hp(\y \C y) = \hp(\x, \z \C x) J_{\x, \z}(\y), \label{trans_2}
\end{eqnarray}
where $X$ is a transformation of $Y$ defined in (\ref{eq:transform}).
Note that the two predictive density functions on the left and right sides of the above  equation may have different functional forms; however, to simplify the
notation, we use the same symbol $\hp$ to represent them when the
context is clear.

Now, the average KL risk can be represented as
\begin{eqnarray}\label{teo2.1:eq1}
R(f, \hp) &=& \EE_{Y, \Y | f} \log \frac{p(\Y \C f)}{\hp(\Y \C Y)} \nonumber \\[-8pt]\\[-8pt]
&=& \EE_{X, \X, \Z | \th} \log \frac{p(\X, \Z \C \th)}{\hp(\X, \Z \C
X)},\nonumber
\end{eqnarray}
where the second equality follows from (\ref{trans_1}) and
(\ref{trans_2}).  Since $\X$ and $\Z$  are independent, we can split $p(\x,\z \C \th)$ as
\begin{equation} \label{teo2.1:eq2}
p(\x,\z \C \th) = p(\x \C \th) p(\z),
\end{equation}
where $p(\z)$ has a known distribution $N_{m-n}(\mathbf{0},
I_{m-n})$ Moreover, to evaluate the minimax risk, it suffices to consider
predictive density estimators in the  form
\begin{equation} \label{teo2.1:eq3}
 \hp(\x,\z \C x) = \hp(\x \C x) p(\z)
\end{equation}
because any predictive density  $\hp(\x,\z \C x)$ can be written as
$\hp(\x,\z \C x) = \hp(\x \C x) \hp(\z \C x, \x)$, and if $\hp(\z \C
x, \x)$ is equal to $p(\z)$, then this density estimator is
dominated by $\hp(\x \C x)p(\z)$, due to the non-negativity of KL
divergence.

Combining (\ref{teo2.1:eq1})--(\ref{teo2.1:eq3}), we have
\[ R(f, \hp) =  \EE_{X, \X | \th} \log \frac{p(\X \C \th)}{\hp(\X; X)} = R(\th, \hp). \]
Consequently, the minimax risk in the non-parametric regression model is equal to the minimax risk in the Gaussian sequence model.
\end{pf*}


\begin{pf*}{Proof of Lemma \ref{lem4.1}}
 Let $\mathcal{Q}$ be the collection of all (generalized) Bayes
 predictive densities.   Then, by
   \cite{r21}, Theorem 5, $\mathcal{Q}$ is a
complete
 class for the problem of predictive density  estimation  under KL
loss.  Therefore, the minimax risk among all possible density
estimators is equivalent to the minimax risk among (generalized)
Bayes estimators, namely,
\begin{eqnarray*}
R(\Theta) = \inf_{\hp} \sup_{\th \in \Th} R(\th, \hp) = \inf_{\hp
\in \mathcal{Q}} \sup_{\th \in \Th} R(\th, \hp).
\end{eqnarray*}

Consider a Gaussian distribution $\pi_S =  N(0, S)$, where $S =
\operatorname{diag}(s_1^2, \dots, s_n^2)$ and the $s_i$'s satisfy condition
(\ref{cond:s_i}). Then,
\begin{eqnarray}
R(\Theta) &=& \inf_{\hp \in \mathcal{Q}} \sup_{\th \in \Th} R(\th, \hp) \\
&\ge& \inf_{\hp \in \mathcal{Q}} \int_{\Th} R(\th, \hp) \pi_S(\th)\, \mathrm{d} \th \nonumber \\
                                                                      &\ge& \inf_{\hp \in \mathcal{Q}} \int_{\RR^n} R(\th, \hp) \pi_S(\th) \,\mathrm{d} \th - \sup_{\hp \in \mathcal{Q}} \int_{\Theta^c} R(\th, \hp) \pi_S(\th)\, \mathrm{d} \th \nonumber \\
                                                                      &\ge& \inf_{\hp \in \mathcal{Q}} \int_{\RR^n} R(\th, \hp) \pi_S(\th) \,\mathrm{d} \th - \sup_{\hp \in \mathcal{Q}} \int_{\Theta^c} R(\th, \hp) \pi_S(\th) \,\mathrm{d} \th. \label{ineq:bound_1}
\end{eqnarray}

The~first term of (\ref{ineq:bound_1}) is the Bayes risk under
$\pi_S$ over the unconstrained parameter space $\RR^n$. It is
achieved by the linear predictive density $\hp_S$; see \cite{r22}.  Therefore,
\begin{eqnarray} \label{ineq:term_1}
\inf_{\hp \in \mathcal{Q}} \int_{\RR^n} R(\th, \hp) \pi_S(\th)\, \mathrm{d} \th &=& \int_{\RR^n} R(\th, \hp_S) \pi_S(\th)\, \mathrm{d} \th \nonumber \\[-8pt]\\[-8pt]
                                                     &=& \frac{n}{2m} \log \frac{\s2_n}{\s2_{n+m}} + \frac{1}{2m}\sum_{i=1}^n \log \frac{\s2_{n+m} + s_i^2}{\s2_n +
                                                     s_i^2}.\nonumber
\end{eqnarray}

To bound the second term of (\ref{ineq:bound_1}), note that for any
Bayes predictive density $\hp_\pi \in \mathcal{Q}$,
\begin{eqnarray}
R(\th, \hp_{\pi}) & = & \frac{1}{m} \EE_{X, \X \mid  \th}  \log \frac{p(\X \vert  \th)}{\int_{\Th} p(\X \vert \th') \pi(\th' \vert  X) \,\mathrm{d} \th'} \nonumber \\
                  &\le& \frac{1}{m} \EE_{X, \X \mid \th} \int_{\Th}\log \frac{p(\X \vert  \th)}{p(\X \vert  \th')} \pi(\th' \vert  X)\, \mathrm{d} \th' \label{ineq:risk:bayes:1}\\
                  & = & \frac{1}{m} \EE_{X \mid \th} \int_{\Th} \frac{\| \th - \th' \|^2}{2 \s2_m} \pi(\th' \vert  X) \, \mathrm{d} \th' \nonumber \\
                  &\le& \frac{1}{m \s2_m} \EE_{X \mid \th} \int_{\Th}  ( \| \th \|^2 + \|\th'\|^2  ) \pi(\th' \vert  X) \,\mathrm{d} \th'  \label{ineq:risk:bayes:2} \\
                  &\le& \frac{1}{m \s2_m}  \biggl(  \| \th \|^2  + \frac{C}{a_1^2} \biggr), \label{ineq:risk:bayes:3}
\end{eqnarray}
where (\ref{ineq:risk:bayes:1}) is due to  Jensen's inequality,
(\ref{ineq:risk:bayes:2}) is due to $\| \th - \th'\|^2 \le 2 \| \th
\|^2 + 2\|\th'\|^2$ and  (\ref{ineq:risk:bayes:3}) is due to
\begin{eqnarray*}
\int_{\Th} \|\th'\|^2 \pi(\th' \C x) \,\mathrm{d} \th' &\le& \sup_{\th' \in \Th}
\| \th' \|^2\\
&  \le& \frac{1}{a_1^2} \sup_{\th \in \Th} \sum_{i=1}^n
a_i^2 \th_i^{'2} = \frac{C}{a^2_1}.
\end{eqnarray*}
Therefore,
\begin{equation} \label{eq:bound:term2}
\sup_{\hp \in \mathcal{Q}} \int_{\Theta^c} R(\th, \hp) \pi_S(\th)\, \mathrm{d}
\th \le \frac{1}{m v_m} \biggl[ \int_{\Th^c} \|\th\|^2 \pi_S(\th)\, \mathrm{d}
\th + \frac{C}{a_1^2} \pi_S(\Th^c) \biggr],
\end{equation}
where $\pi_S(\Theta^c) = \int_{\Th^c} \pi_S(\th)\, \mathrm{d} \th$.  Using
the Cauchy--Schwarz inequality, we can further bound the right-hand side of
(\ref{eq:bound:term2}) as follows:
\begin{eqnarray*}
&& \frac{1}{m v_m} \biggl[ \int_{\Th^c} \|\th\|^2 \pi_S(\th) \,\mathrm{d} \th + \frac{C}{a_1^2} \pi_S(\Th^c) \biggr]  \\
&&\quad \le \frac{1}{m v_m} \Biggl[ \sum_{i=1}^n \biggl( \int_{\Th^c} \th_i^4 \pi_S(\th) \,\mathrm{d} \th  \biggr)^{1/2} \sqrt{\pi_S(\Th^c)} + \frac{C}{a_1^2} \pi_S(\Th^c) \Biggr]  \\
                                                                        &&\quad  =  \frac{1}{m v_m} \Biggl[ \sqrt{3} \sqrt{\pi_S(\Th^c)} \sum_{i=1}^n s_i^2 +  \frac{C}{a_1^2} \pi_S(\Th^c) \Biggr] \\
                                                                        &&\quad \le \frac{1}{m v_m} \biggl[ \sqrt{3}\frac{C}{a_1} \sqrt{\pi_S(\Th^c)} +
                                                                        \frac{C}{a_1} \pi_S(\Th^c) \biggr].
\end{eqnarray*}
Then, by
\cite{r2}, Proposition 2, which states that if $\epsilon_1, \dots, \epsilon_m$ are independent
Gaussian random variables with $E \epsilon_k = 0$ and $E
\epsilon^2_k = \sigma_k^2$, then
\[ \mathbb{P}  \Biggl( \sum_{k=1}^m \epsilon_k^2 > Q  \Biggr) \le \exp
   \biggl\{ - \frac{(Q - \sum_{k=1}^m \sigma^2_k)^2}{4 \sum_{k=1}^m
  \sigma_k^4}  \biggr\},  \]
we have
\begin{eqnarray} \label{order}
\sqrt{\pi_S(\Th^c)} =  \Biggl[ \mathbb{P}  \Biggl( \sum_{i=1}^n a_i^2
\th_i^2 > C  \Biggr)  \Biggr]^{1/2} \le \s2_n^{ \alpha},
\end{eqnarray}
due to condition (\ref{cond:s_i}).

Combining (\ref{ineq:bound_1}), (\ref{ineq:term_1}), (\ref{eq:bound:term2}) and (\ref{order}), the theorem then follows immediately.
\end{pf*}
\end{appendix}

\section*{Acknowledgements}
The~authors would like to thank Edward I. George for helpful
discussions and the Associate Editor for generous insights and suggestions.  This work was supported in part by the
National Science Foundation under award numbers DMS-07-32276 and DMS-09-07070.  Any opinions, findings and conclusions
or recommendations expressed in this material are those of the authors and do not necessarily reflect the views of
the National Science Foundation.

\printhistory


\begin{thebibliography}{99}
\bibitem{r22}
{Aitchison, J.} (1975).  Goodness of prediction fit.  \textit{Biometrika} \textbf{62} 547--554.

\bibitem{r1}
{Aslan, M.} (2006). Asymptotically minimax Bayes predictive
densities. \textit{Ann. Statist.} \textbf{34} 2921--2938.

\bibitem{r2}
{Belitser, E.N.} and {Levit, B.Y.} (1995). On minimax
filtering over ellipsoids. \textit{Math. Methods
Statist.} \textbf{3} 259--273.

\bibitem{r3}
{Belitser, E.N.} and {Levit, B.Y. } (1996).
Asymptotically minimax nonparametric regression in $L_2$.
\textit{Statistics} \textbf{28} 105--122.
\MR{1405604}

\bibitem{r21} {Brown, L.D.}, {George, E.I.} and
{Xu, X.} (2008). Admissible predictive density estimation.
\textit{Ann. Statist.} \textbf{36} 1156--1170.
\MR{2418653}

\bibitem{r4}
{Brown, L.D.} and {Low, M.G.} (1996). Asymptotic
equivalence of nonparametric regression and white noise. \textit{Ann. Statist.} \textbf{24} 2384--2398.
\MR{1425958}

\bibitem{r5}
{Diaconis, P.} and {Ylvisaker, D.} (1979). Conjugate
priors for exponential families.  \textit{Ann. Statist.} \textbf{7} 269--281.
\MR{0520238}

\bibitem{r6}
{Efromovich, S.Y.} (1994). On adaptive estimation of
nonlinear functionals.  \textit{Statist. Probab. Lett.} \textbf{19}
57--63.
\MR{1253313}

\bibitem{r7}
{Efromovich, S.Y.} (1999). \textit{Nonparametric Curve
Estimation: Methods, Theory and Applications}. New York:  Springer.
\MR{1705298}

\bibitem{r8}
{Efromovich, S.Y.} and {Pinsker, M.S.} (1982).
Estimation of square-integrable probability density of a random
variable.   \textit{Probl. Inf. Transm.} \textbf{18} 175--189.
\MR{0711898}

\bibitem{r9}
{George, E.I.}, {Liang, F.} and {Xu, X.}
(2006). Improved minimax prediction under Kullback--Leibler loss.
\textit{Ann. Statist.} \textbf{34} 78--91.
\MR{2275235}

\bibitem{r10}
{George, E.I.} and {Xu, X.} (2008). Predictive density
estimation for multiple regression. \textit{Econometric Theory} \textbf{24} 1--17.
\MR{2391619}

\bibitem{r11}
{Goldenshluger, A.} and {Tsybakov, A.B.} (2003).
Optimal prediction for linear regression with infinitely many
parameters. \textit{J. Multivariate Anal.} \textbf{84}
40--60.
\MR{1965822}

\bibitem{r12}
{Golubev, G.K.} and {Nussbaum, M.} (1990). A~risk
bound in Sobolev class regression. \textit{Ann. Statist.} \textbf{18} 758--778.
\MR{1056335}

\bibitem{r13}
{Ibragimov, I.A.} and {Has'minskii, R.Z.} (1977). On
the estimation of an infinite-dimensional parameter in Gaussian
white noise.  \textit{Soviet Math. Dokl.} \textbf{236} 1053--1055.
\MR{0483232}

\bibitem{r14}
{Ibragimov, I.A.} and {Has'minskii, R.Z.} (1984). On
nonparametric estimation of values of a linear functional in a
Gaussian white noise.  \textit{Teor. Veroyatn. Primen.} \textbf{29} 19--32.
\MR{0739497}

\bibitem{r15}
{Johnstone, I.M.} (2003). Function estimation and
Gaussian sequence models.  Draft of a Monograph.

\bibitem{r16}
{Liang, F.} and {Barron, A.} (2004). Exact minimax
strategies for predictive density estimation, data compression and
model selection. \textit{IEEE Trans. Inform. Theory} \textbf{50}
2708--2726.
\MR{2096988}

\bibitem{r17}
{Nussbaum, M.} (1996). Asymptotic equivalence of density
estimation and Gaussian white noise. \textit{Ann. Statist.} \textbf{24} 2399--2430.
\MR{1425959}

\bibitem{r18}
{Nussbaum, M.} (1999). Minimax risk: Pinsker bound. In
\textit{Encyclopedia of Statistical Sciences} (S. Kotz, ed.) 451--460.
New York: Wiley.

\bibitem{r19}
{Pinsker, M.S.} (1980). Optimal filtering of square
integrable signals in Gaussian white noise. \textit{Probl. Inf.
Transm.} \textbf{2} 120--133.

\bibitem{r20}
{Tsybakov, A.B.} (1997). On nonparametric estimation of
density level sets.  \textit{Ann. Statist.} \textbf{25}
948--969.
\MR{1447735}

\end{thebibliography}
\end{document}